\documentclass[11pt]{article}
\usepackage[active]{srcltx}
\usepackage{amssymb,amsmath}
\usepackage{graphicx}

\newtheorem{proposition}{Proposition}[section]
\newtheorem{theorem}{Theorem}[section]

\newtheorem{definition}{Definition}[section]
\newtheorem{remark}{Remark}[section]

\numberwithin{equation}{section}

\begin{document}

\begin{center}{\large\sc
Metric Projections  versus Non-Positive Curvature%. Part I. Theory
%: from CAT(0)-spaces to Finsler manifolds
%Metric Projections.
%Part I.\\ Old and New Results in Riemann-Finsler Geometry
}\\
\vspace{0.5cm}
%by \vspace{0.5cm}

{\large  Alexandru Krist\'aly%\footnote{Research supported by grant
%CNCSIS PN-II-ID-PCE-2011-3-0241.}
\\
  {\normalsize Department of Economics, Babe\c s-Bolyai
University,  400591 Cluj-Napoca, Romania
%\\{\footnotesize Email address: alexandrukristaly@yahoo.com}
}\\
\vspace{0.5cm}
 Du\v san Repov\v s\\  {\normalsize Faculty  of Education, and Faculty of Mathematics and Physics, University of
 Ljubljana, P.O.B. 2964, 1001
 Ljubljana,  Slovenia
 %\\
 %{\footnotesize Email address: dusan.repovs@guest.arnes.si}
 }\\
   }
\end{center}

%\vspace{1cm}

\begin{abstract}
{\footnotesize \noindent In this paper two metric properties on
geodesic length spaces are introduced by means of the metric
projection, studying their validity on Alexandrov and Busemann NPC
spaces. In particular, we prove that both properties characterize
the non-positivity of the sectional curvature on Riemannian
manifolds. Further results are also established on
reversible/non-reversible Finsler-Minkowski spaces.

%  We show that
%non-positively curved Alexandrov spaces
% verify the {\it double-pro\-jec\-tion
%property}, i.e., a point is a best approximation element between two
%'small' geodesic convex, compact sets $S_1$ and $S_2$ if and only if
%it is a fixed point of the metric projection map $P_{S_1}\circ
%P_{S_2}$. Although the converse need not hold in generic geodesic
%length spaces (e.g., on reversible Finslerian Minkowski spaces), the
%local double-projection property is equivalent to the non-positivity
%of the sectional curvature on $n$-dimensional Riemannian manifolds,
%$n\geq 2$.
%The tools
%used in the proofs are the fundamental inequality for Minkowski
%norms known from Finsler geometry and the Ptolemaic inequality on
%CAT(0)-spaces. A curvature rigidity result in the Riemannian setting
%is also provided.
}
%A  metric space $(M,d)$ verifies the {\it local double-projection
%property} if every point $p\in M$ has a neighborhood $U\subset M$
%such that $(U,d)$ is a geodesic length space, and for every two
%geodesic convex, closed, $U$-Chebisev sets $S_1,S_2\subset U$,  the
%element $q\in S_1$ is the best approximation point from $S_1$ to
%$S_2$ if and only if $q$ is a fixed point of $P_{S_1}\circ P_{S_2}$,
%where $P_S$ is the metric projection to the set $S$.
\end{abstract}


\section{Introduction}\label{sect:1}

The curvature notions on geodesic length spaces are formulated in
terms of the metric distance. Most of them refer to non-positively
curved spaces (shortly, NPC spaces) defined by means of certain
metric inequalities. Here, we recall (non-rigorously) three such
notions:

\begin{itemize}
  \item[(a)] {\it Alexandrov NPC spaces} (see \cite{Alexandrov}):
  small geodesic triangles are thinner than their Euclidean
  comparison triangles;
  \item[(b)] {\it Busemann NPC spaces} (see \cite{Busemann}): in
  small geodesic triangles the geodesic segment connecting the
  midpoints of two sides is at most half as long as the third
  side;
  \item[(c)] {\it Pedersen NPC spaces} (see \cite{Pedersen}): small capsules (i.e.,
the loci equidistant to geodesic segments) are geodesic convex.
\end{itemize}

\noindent It is well-known that

"Alexandrov NPC spaces $\subset$ Busemann NPC spaces $\subset$
Pedersen NPC spaces,"

\noindent where the inclusions are proper in general. However, on
Riemannian manifolds, all these curvature notions coincide, and they
characterize the non-positivity of the sectional curvature. For
systematic presentation of NPC spaces, we refer the reader to the
monographs of Bridson and Haefliger \cite{Br-Ha}, Busemann
\cite{Busemann}, and Jost \cite{Jost}.

The aim of our paper is to capture new features of non-positively
curved geodesic length spaces by means of the {\it metric
projection} map. Roughly speaking, on a metric space $(M,d)$, we
consider the following two properties we are dealing with in the
sequel (for precise notions, see Definitions \ref{def-double-proj}
\& \ref{def-non-exp}):
\begin{itemize}
  \item[(I)] {\it Double-pro\-jec\-tion
property}: a point is the best approximation element between two
small geodesic convex sets $S_1,S_2\subset M$ if and only if it is a
fixed point of the metric projection map $P_{S_1}\circ P_{S_2}$.
  \item[(II)] {\it Projection non-expansiveness property}: the metric projection map $P_S$ is
  non-expansive for small
  geodesic convex sets $S\subset M$.
\end{itemize}
Our results can be summarized as follows (for precise statements and
detailed comments, see Section \ref{sect:22}). Although Busemann NPC
spaces do not satisfy in general the above properties (see Remark
\ref{Busemann-remark}), Alexandrov NPC spaces satisfy both of them
(see Theorem \ref{Ptolemaic-theorem}, and Bridson and Haefliger
\cite[Proposition 2.4]{Br-Ha}). Furthermore, generic
Finsler-Minkowski spaces satisfy the global double-projection
property (see Theorem \ref{Minkowski-theorem}), but not the global
projection non-expansiveness property. Finally, we prove that both
properties (I) and (II) encapsulate the concept of non-positive
curvature in the Riemannian context; namely, we prove that for
Riemannian manifolds the double-projection property, the projection
non-expansiveness property and the non-positivity of the sectional
curvature are equivalent conditions (see Theorem
\ref{Riemann-double-proj}).

\section{Main results and remarks}\label{sect:22}

Let $(M,d)$ be a metric space and let
\begin{equation}\label{projection-definition}
    P_S(q)=\{s\in S:d(q,s)=\inf_{z\in S}d(q,z)\}
\end{equation}
be the usual {\it metric projection} of the point $q\in M$ onto the
nonempty set $S\subset M.$ If $S\subset U\subset M$, the set $S$ is
called $U-${\it proximinal} if $P_S(q)\neq \emptyset$  for every
$q\in U$ (w.r.t. the metric $d$), and $U-${\it Chebishev} if
$P_S(q)$ is a singleton for every $q\in U.$

\begin{definition}\label{def-double-proj}
The metric space $(M,d)$ satisfies the {\rm double-pro\-jec\-tion
property} if every point $p\in M$ has a neighborhood $U\subset M$
such that $(U,d)$ is a geodesic length space, and for every two
geodesic convex, $U-$proximinal
%$U-$Chebisev
sets $S_1,S_2\subset U$ and for some $q\in S_1$ the following
statements are equivalent:

$(DP_1)$ $q\in  (P_{S_1}\circ P_{S_2})(q)$;

$(DP_2)$ there exists $\tilde q\in P_{S_2}(q)$ such that $d(q,\tilde
q)\leq d(z_1,z_2)$ for all   $z_1\in
  S_1$, $z_2\in S_2.$

\noindent If $U=M$, then $(M,d)$ satisfies the {\rm global
double-pro\-jec\-tion property.}
\end{definition}

\noindent The element $q\in S_1$ satisfying $(DP_2)$ is called the
{\it best approximation point from the set $S_1$ to $S_2$}. We
notice that $P_{S_i}$ ($i=1,2$) may be set-valued maps in the
Definition \ref{def-double-proj}.

\begin{definition}\label{def-non-exp}
The metric space $(M,d)$ satisfies the {\rm projection
non-expansiveness property} if every point $p\in M$ has a
neighborhood $U\subset M$ such that $(U,d)$ is a geodesic length
space, and for every geodesic convex, $U-$proximinal set $S\subset
U$,  one has
\begin{equation}\label{non-exp}
    d(P_S(q_1),P_S(q_2))\leq d(q_1,q_2)\ {for\ every}\ q_1,q_2\in U.
\end{equation}
If $U=M$, then $(M,d)$ satisfies the {\rm global projection
non-expansiveness property.}
\end{definition}

\noindent Note that if a set $S$ satisfies (\ref{non-exp}), it is
necessarily a $U-$Chebishev set.

\begin{remark}\label{Busemann-remark}\rm
Let us discuss first the relationship between these properties and
Busemann NPC spaces. We recall that every Minkowski space in the
classical sense (i.e., normed linear space with strictly convex unit
ball) is a Busemann NPC space, see Busemann \cite{Busemann}.

(a) {\it Double-projection property fails in Busemann NPC spaces}:
Let $(\mathbb R^3,F)$ be a Min\-kow\-ski space with strictly convex
unit balls. Assume that $F$ is non-differentiable at $p\in I=\{p\in
\mathbb R^3:F(p)=1\}$; then due to the symmetry of $F$, the same
holds at $q=-p\in I$. On account of this assumption, we may consider
supporting planes $H_{p}$ and $H_q$ at $p$ and $q$ to the unit ball
$B=\{p\in \mathbb R^3:F(p)\leq 1\}$, respectively, such that
$H_p\cap H_q\neq \emptyset$. We translate $H_q$ to the origin,
denoting it by $H_0$. Let us finally consider an arbitrary plane $H$
containing the origin and the point $p$, and  $H_0\cap H_p\cap
H=\{z\}$. If $S_1=[p,z]$ and $S_2=[0,z]$, then by construction, one
has $P_{S_1}(0)=p$ and $P_{S_2}(p)=0,$ thus $(P_{S_2}\circ
P_{S_1})(0)=0.$ If the double-projection property holds (up to a
scaling of the indicatrix $I$), then we have that $d_F(0,p)\leq
d_F(z_1,z_2)$ for every $z_1\in S_1$ and $z_2\in S_2$. Let
$z_1=z_2=z\in S_1\cap S_2$; the latter inequality implies the
contradiction $1=F(p)=d_F(0,p)\leq 0$.

(b) {\it Global projection non-expansiveness property fails in
Busemann NPC spaces}: Due to Phelps \cite[Theorem 5.2]{Phelps}, a
Minkowski space (with dimension at least three) which satisfies the
global projection non-expansiveness property, is necessarily
Euclidean.
\end{remark}

Next, we treat these two properties on Alexandrov NPC spaces. First,
it is a well known fact that every Alexandrov NPC space satisfies
the projection non-expansiveness property, see Bridson and Haefliger
\cite[Proposition 2.4]{Br-Ha}. Our first result reads as follows.

%Standard arguments show that the (local) double-projection property
%holds in the Euclidean case.
 %In general, the following result
%can be stated.

\begin{theorem}\label{Ptolemaic-theorem}
Every Alexandrov NPC space satisfies the  double-projection
property.
\end{theorem}

\begin{remark}\rm\label{remark-bizonyitasi-modok}
(a) We provide two independent proofs of Theorem
\ref{Ptolemaic-theorem}, each of them exploiting basic properties of
Alexandrov NPC spaces: (1) Pythagorean and Ptolemaic inequalities;
(2) the first variation formula and non-expansiveness of the metric
projection.

(b) With respect to Remark \ref{Busemann-remark}, if we assume that
a Busemann NPC space is also Ptolemy (i.e., the Ptolemaic inequality
holds for every quadruple), the double-projection property holds. In
fact, the latter statement is precisely Theorem
\ref{Ptolemaic-theorem}, exploiting the famous characterization of
CAT$(0)-$spaces by Foertsch, Lytchak and Schroeder \cite{FLS}, i.e.,
a metric space is a CAT$(0)-$space if and only if it is a Ptolemy
and a Busemann NPC space.
\end{remark}
We now present a genuinely different class of spaces where the
double-projection property holds.

\begin{theorem}\label{Minkowski-theorem}
Every reversible Finsler-Minkowski space satisfies the global
double-pro\-jec\-tion property.
\end{theorem}

\begin{remark}\rm (a) Hereafter, the Finsler-Minkowski space is understood in
the sense of Finsler geometry, see Bao, Chern and Shen \cite{BCS};
in particular, we assume that the norm $F$ belongs to $ C^2(\mathbb
R^n\setminus \{0\});$ see Section \ref{sect:2}. As we already
pointed out in Remark \ref{Busemann-remark}(a), the
double-projection fails on Minkowski spaces with non-differentiable
unit balls.

(b) We emphasize that the proof of Theorem \ref{Minkowski-theorem}
cannot follow any of the lines described in Remark
\ref{remark-bizonyitasi-modok}(a). First, a rigidity result due to
Schoenberg \cite{Schoenberg} shows that any Minkowski space on which
the Ptolemaic inequality holds is necessarily Euclidean; see also
Buckley, Falk and Wraith \cite{BFW}. Second, if we want to apply the
global projection non-expansiveness property, we come up against the
rigidity result of Phelps \cite[Theorem 5.2]{Phelps}, see also
Remark \ref{Busemann-remark}(b).
 However, the
fundamental inequality of Finsler geometry and some results from
Krist\'aly, R\u adulescu and Varga \cite{KRV} provide a simple proof
of Theorem \ref{Minkowski-theorem}, where the fact that $F$ belongs
to $C^2(\mathbb R^n\setminus \{0\})$ plays an indispensable role.
\end{remark}

In spite of the above remarks, the following characterization can be
proved in the Riemannian framework which entitles us to assert that
the notions introduced in Definitions \ref{def-double-proj} \&
\ref{def-non-exp} provide new features of the non-positive
curvature.

\begin{theorem}\label{Riemann-double-proj}
Let $(M,g)$ be a smooth Riemannian manifold and $d_g$ the induced
metric on $M$. Then the following assertions are equivalent:
\begin{itemize}
  \item[{\rm (i)}] $(M,d_g)$ satisfies the  double-projection
  property;
  \item[{\rm (ii)}] $(M,d_g)$ satisfies the projection non-expansiveness
  property;
  \item[{\rm (iii)}] the sectional curvature of $(M,g)$ is non-positive.
\end{itemize}
\end{theorem}
The proof of Theorem \ref{Riemann-double-proj} is based on the
Toponogov comparison theorem and on the formula of the sectional
curvature given by the Levi-Civita parallelogramoid.

 In order for the paper to be
self-contained, we recall in Section \ref{sect:2}  some basic
notions and results from  Alexandrov NPC spaces and
Finsler-Minkowski spaces. In Section \ref{sect-proof} we present the
proof of Theorems \ref{Ptolemaic-theorem} and
\ref{Riemann-double-proj}, while in Section \ref{sect-4} we prove
Theorem \ref{Minkowski-theorem} and also discuss  some aspects of
the double-projection property on non-reversible Finsler-Minkowski
spaces.

\section{Preliminaries%: Alexandrov NPC spaces and Min\-kowski spaces
}\label{sect:2}

\noindent {\bf A. Alexandrov NPC spaces.} We recall those notions
and results from the theory Alexandrov NPC spaces which will be used
in the proof of Theorems \ref{Ptolemaic-theorem} and
\ref{Riemann-double-proj}; for details, see Bridson and Haefliger
\cite[Chapter II]{Br-Ha}, and Jost \cite{Jost}.

A metric space $(M,d)$ is a {\it geodesic length space} if for every
two points $p,q\in M$, there exists the shortest geodesic segment
joining them, i.e., a continuous curve $\gamma:[0,1]\to M$ with
$\gamma(0)=p$, $\gamma(1)=q$ and $l(\gamma)=d(p,q)$, where
$$l(\gamma)=\sup\left\{\sum_{i=1}^m d(\gamma(t_{i-1}),\gamma(t_{i})):0=t_0<...<t_m=1, \ m\in \mathbb N \right\}.$$ We assume that geodesics are parametrized proportionally
by the arc-length.

Given a real number $\kappa$, let $M_\kappa^2$ be the
two-dimensional space form with curvature $\kappa,$ i.e.,
$M_0^2=\mathbb R^2$ is the Euclidean plane, $M_\kappa^2$ is the
sphere with radius $1/\sqrt{\kappa}$ if $\kappa>0$, and $M_\kappa^2$
is the hyperbolic plane with the function multiplied by
$1/\sqrt{-\kappa}$ if $\kappa <0$. If $p,q,r\in M$, a {\it geodesic
triangle} $\Delta(p,q,r)$ in $(M,d)$ is defined by the three
vertices and a choice of three sides which are geodesic segments
joining them (they need not be unique). A triangle $\overline
\Delta(\overline p,\overline q,\overline r)\subset M_\kappa^2$ is a
{\it comparison triangle} for $\Delta(p,q,r)\subset M$, if
$d(p,q)=d(\overline p,\overline q)$, $d(p,r)=d(\overline p,\overline
r)$, and $d(r,q)=d(\overline r,\overline q)$. If
$d(p,q)+d(q,r)+d(r,p)<2 D_\kappa$ (where $D_\kappa={\rm
diam}(M_\kappa^2)$), such a comparison triangle exists and it is
unique up to isometries.   A point $\overline x \in {\rm
Im}(\overline \gamma)$ is a comparison point for $ x \in {\rm
Im}(\gamma)$ if $d(p,x)=d(\overline p,\overline x)$, where
$\gamma:[0,1]\to M$ and $\overline \gamma:[0,1]\to M_\kappa^2$ are
geodesic segments such that $\gamma(0)=p$, $\gamma(0)=\overline p$,
and $l(\gamma)=l(\overline \gamma).$

Let $\Delta(p,q,r)\subset M$ be a geodesic triangle with perimeter
less than $2D_\kappa$,  and let $\overline \Delta(\overline
p,\overline q,\overline r)\subset M_\kappa^2$ be its comparison
triangle. The triangle $\Delta(p,q,r)$ satisfies the {\it
CAT$(\kappa)-$in\-equality}, if for every $x,y\in
\Delta(p,q,r)$, %=\gamma_{p,q}\cup \gamma_{p,r}\cup \gamma_{r,q}$ (where
%$\gamma_{p,q}$ denotes the geodesic segment joining $p$ and $q$,
%etc.),
for the comparison points $\overline x,\overline y\in \overline
\Delta(\overline p,\overline q,\overline r)$ one has $d( x,y)\leq
d(\overline x,\overline y).$ The geodesic length space $(M,d)$ is a
{\it CAT$(\kappa)-$space} if all geodesic triangles in $M$ with
perimeter less than $2D_\kappa$ satisfy the
CAT$(\kappa)-$inequality. The metric space $(M,d)$ is an {\it
Alexandrov NPC space} if it is locally a CAT$(0)-$space, i.e., for
every $p\in M$ there exists $\rho_p>0$ such that $B(p,\rho_p)=\{q\in
M:d(p,q)<\rho_p\}$ is a CAT$(0)-$space.

 A set $S\subset M$ is {\it geodesic convex} if for every
two points $p,q\in S$, there exists a unique geodesic segment
joining $p$ to $q$ whose image is contained in $S$. The projection
map $P_S:M\to 2^S$ is defined by (\ref{projection-definition}).

%We collect some The following properties of CAT(0)-spaces are
%well-known.

\begin{proposition}\label{CAT0-prop}
Let $(M,d)$ be a {\rm CAT$(0)$}$-$space.  Then the following
properties hold:

\begin{itemize}
\item[{\rm (i)}] {\rm (See \cite[Proposition 2.2]{Br-Ha})} The
  distance function $d$ is convex.

   \item[{\rm (ii)}] {\rm Projections} {\rm (See \cite[Proposition 2.4]{Br-Ha}):}  If $S\subset M$ is a geodesic
  convex $M-$proximinal  set,   then it is
  $M-$Chebishev, i.e., $P_S(q)$ is a singleton for every $q\in M.$
  Moreover, $P_S$ is non-expansive, i.e., {\rm (\ref{non-exp})}
  holds  on $M.$ If $q\notin S$ and $z\in S$, then $\angle_{P_S(q)}(q,z)\geq \pi/2,$
  where $\angle_p(z_1,z_2)$ denotes the Alexandrov angle between the
  unique geodesic segments joining $p$ to $z_1$ and $z_2,$
  respectively.

  \item[{\rm (iii)}] {\rm First variation formula} {\rm (See \cite[Corollary 3.6]{Br-Ha}):}  If $\gamma:[0,1]\to M$
  is a
  geodesic segment with $\gamma(0)=p$, and $z\in M$ is a distinct point from $p$, then
  $$\cos \angle_p(\gamma(t),z)=\lim_{s\to 0^+}\frac{d(p,z)-d(\gamma(s),z)}{s},\ t\in (0,1].$$

  \item[{\rm (iv)}] {\rm Pythagorean inequality (See \cite[Theorem 2.3.3]{Jost}):} If $p\in M,$ $\gamma:[0,1]\to M$
  is a geodesic segment,  and $\gamma(0)=P_{{\rm
  Im}(\gamma)}(p)$, then
$$d^2(p,\gamma(0))+d^2(\gamma(0),\gamma(1))\leq d^2(p,\gamma(1)) .$$
\vspace{-0.6cm}
  \item[{\rm (v)}] {\rm Ptolemaic inequality  (See \cite{FLS, Kay}):} For every quadruple $q_i\in M$, $i=1,...,4$,
one has
$$d(q_1,q_3)\cdot d(q_2,q_4)\leq d(q_1,q_2)\cdot
d(q_3,q_4)+d(q_1,q_4)\cdot d(q_2,q_3).$$

\end{itemize}
\end{proposition}

\begin{remark}\rm  If $P_S(q)$ is a singleton for some $q\in M$, we do
not distiguish between the set and its unique point.
\end{remark}

\noindent  {\bf B. Finsler-Minkowski spaces.} Let $F:\mathbb R^n\to
[0,\infty)$ be a {\it positively} {\it homogenous Minkowski norm},
i.e., $F$ satisfies the properties:

(a) $F\in C^2(\mathbb R^n\setminus \{0\});$

(b) $F(ty)=tF(y)$ for all $t\geq 0$ and $y\in \mathbb R^n;$

(c) The Hessian matrix %$h_{ij}(\xi)= [\frac{1}{2}H^2]_{\xi_i\xi_j}(\xi)$
  $g_y=\nabla^2 ({F^2}/{2})(y)$ is positive definite for all $y\neq 0.$

\noindent The Minkowski norm  $F$ is said to be {\it absolutely
homogeneous} if in addition, we have

(b') $F(ty)=|t|F(y)$ for all $t\in \mathbb R$ and $y\in \mathbb
R^n$.

\noindent If (a)-(c) hold, the pair $(\mathbb R^n,F)$ is a {\it
Finsler-Minkowski space}, see Bao, Chern and Shen \cite[\S
1.2]{BCS}, which is the simplest (not necessarily reversible)
geodesically complete Finsler manifold whose flag curvature is
identically zero, the geodesics are straight lines, and the
intrinsic distance between two points $p,q\in \mathbb R^n$ is given
by
\begin{equation}\label{F-metrika}
    d_F(p,q)=F(q-p).
\end{equation}
In fact, $(\mathbb R^n,d_F)$ is a quasi-metric space and in general,
$d_F(p,q)\neq d_F(q,p)$. In particular, $g_{y}=g_{-y}$ for all
$y\neq 0$ if and only if $F$ is absolutely homogeneous; if so,
 $(\mathbb R^n,F)$ is a {\it reversible} Finsler-Minkowski space.

Let $S\subset \mathbb R^n$ be a nonempty set. Since $(\mathbb
R^n,F)$ is not necessarily reversible, we define the {\it forward}
(resp. {\it  backward}) {\it  metric projection}s of $q$ to $S$ as
follows:

\begin{itemize}
\item $P_{S}^{+}(q)=\left\{  s_{f}\in S:d_{F}(q,s_{f})=\inf_{s\in S}
d_{F}(q,s)\right\}  ;$

\item $P_{S}^{-}(q)=\left\{  s_{b}\in S:d_{F}(s_{b},q)=\inf_{s\in S}d_F(s,q)\right\}
$.
\end{itemize}
%A set $S\subset \mathbb R^n$ is {\it convex} if it is convex in the
%usual sense.

\begin{proposition}\label{Minkowski-prop}
Let $(\mathbb R^n,F)$ be a $($not necessarily reversible$)$
Finsler-Minkowski space. Then the following properties hold:

\begin{itemize}
  \item[{\rm (i)}] {\rm (See \cite[Theorem 15.8]{KRV})} If $S\subset \mathbb R^n$ is
   convex and $\mathbb R^n-$proximinal, then $S$ is both forward and backward $\mathbb R^n-$Chebishev, i.e.,
 $P_S^+(q)$ and $P_S^-(q)$ are singletons for every $q\in \mathbb R^n.$
  \item[{\rm (ii)}] {\rm (See \cite[Theorem 15.7]{KRV})} If $S\subset \mathbb R^n$ is
  closed and convex, then
\begin{itemize}
\item[$\bullet$] $s\in P_{S}^{+}(q)$ if and only if $g_{s-q}(s-q,z-s)\geq
0$ for all $z\in S;$

\item[$\bullet$] $s\in P_{S}^{-}(q)$ if and only if $g_{q-s}(q-s,z-s)\leq
0$ for all $z\in S.$
\end{itemize}
%\vspace{-0.5cm}
  \item[{\rm (iii)}] {\rm Fundamental
 inequality of
Finsler geometry (See \cite[p. 6-10]{BCS}):}   For every $y\neq0\neq
w,$ one has
$$ |g_y(y,w)|\leq \sqrt{g_y(y,y)}\cdot \sqrt{g_w(w,w)}=F(y)\cdot
  F(w).$$
\end{itemize}
\end{proposition}

%\begin{remark}\rm  If $(\mathbb R^n,F)$ is reversible,
%$P_S^+(q)=P_S^-(q)$ for every $q\in \mathbb R^n$ and
%\end{remark}

% A Riemannian manifold
%is said to be of {\it Hadamard-type}, if it is simply
%connected, complete, and has non-positive sectional curvature.\\

\section{Proof of Theorems
\ref{Ptolemaic-theorem} and
\ref{Riemann-double-proj}}\label{sect-proof}

\noindent {\bf Proof of Theorem \ref{Ptolemaic-theorem}.} Let $p\in
M$ be fixed. Since $(M,d)$ is an Alexandrov NPC space, there exists
$\rho_p>0$ small enough such that $B(p,\rho_p)$ is a CAT$(0)-$space.
We fix arbitrary two geodesic convex $B(p,\rho_p)-$proximinal sets
$S_1, S_2\subset B(p,\rho_p)$. According to Proposition
\ref{CAT0-prop}(ii), $S_1,S_2$ are $B(p,\rho_p)-$Chebishev sets. We
will prove that $(DP_1)$ is equivalent to $(DP_2)$. Let $q\in S_1.$

{\it Step 1.} $"(DP_2)\Rightarrow (DP_1)"$. Let us choose
$z_2=P_{S_2}(q)\in S_2$ in $(DP_2)$. Therefore, it follows that
$d(q,P_{S_2}(q))\leq d(z_1,P_{S_2}(q))$ for all $z_1\in
  S_1$, which implies that $q\in P_{S_1}(P_{S_2}(q))$. Since  $S_1$ is  $B(p,\rho_p)-$Chebishev, the claim follows.

{\it Step 2.} $"(DP_1)\Rightarrow (DP_2)"$.  Since $S_1$ and $S_2$
are $B(p,\rho_p)-$Chebishev sets, we may assume that $(P_{S_1}\circ
P_{S_2})(q)=q$ in $(DP_1)$. Furthermore,  there exists a unique
element $\tilde q\in S_2$ with $P_{S_2}(q)=\tilde q$ and
$P_{S_1}(\tilde q)=q.$ We shall assume that $d(q,\tilde q)>0;$
otherwise, $(DP_2)$ trivially holds. Fix $z_1\in S_1$ and $z_2\in
S_2$ arbitrarily. Applying the Pythagorean inequality (see
Proposition \ref{CAT0-prop}(iv)) to the point $\tilde q$ and the
geodesic segment joining $q$ to $z_1$, we have
\begin{equation}\label{pit-1}
d^2(q,\tilde q)+d^2(q,z_1)\leq  d^2(z_1,\tilde q).
\end{equation}
In a similar way, one has
\begin{equation}\label{pit-2}
d^2(\tilde q, q)+d^2(\tilde q,z_2)\leq d^2(z_2,q).
\end{equation}
Since $(B(p,\rho_p),d)$ is Ptolemaic (see Proposition
\ref{CAT0-prop}(v)), for the quadruple of points $z_1,z_2,\tilde q,
q\in B(p,\rho_p)$, we obtain
\begin{equation}\label{utolso-remelem}
  d(z_1,\tilde q)\cdot d(z_2,q)\leq d(z_1,z_2)\cdot d(\tilde
q,q)+d(z_1,q)\cdot d(z_2,\tilde q).
\end{equation}
 Assume to the  contrary
that $d(z_1,z_2)<d( q,\tilde q)=d(q,P_{S_2}(q)).$ Then, relation
(\ref{utolso-remelem}) yields $$d(z_1,\tilde q)\cdot d(z_2,q)< d^2(
q,\tilde q)+d(z_1,q)\cdot d(z_2,\tilde q).$$ Combining this relation
with (\ref{pit-1}) and (\ref{pit-2}), we obtain $$[d^2(q,\tilde
q)+d^2(q,z_1)]\cdot[d^2(\tilde q, q)+d^2(\tilde q,z_2)]< [d^2(
q,\tilde q)+d(z_1,q)\cdot d(z_2,\tilde q)]^2,$$ which is equivalent
to $[d(q,z_1)-d(\tilde q,z_2)]^2<0,$ a contradiction. Therefore, we
have  $$d(q,P_{S_2}(q))=d( q,\tilde q)\leq d(z_1,z_2),$$ which
concludes the proof. $\hfill $ $\diamondsuit$

\begin{remark}\rm
For $"(DP_1)\Rightarrow (DP_2)"$ we can give an alternative proof.
As above, let  $\tilde q\in S_2$ with $P_{S_2}(q)=\tilde q$ and
$P_{S_1}(\tilde q)=q$, and fix $z_1\in S_1$ and $z_2\in S_2$
arbitrarily. Let $\gamma:[0,1]\to M$ be the unique geodesic joining
$q=\gamma(0)$ and $\tilde q=\gamma(1)$.  We claim that
\begin{equation}\label{propro}
    P_{{\rm Im}(\gamma)}(z_1)=q\ {\rm and}\ P_{{\rm
Im}(\gamma)}(z_2)=\tilde q.
\end{equation}
Since $P_{S_1}(\tilde q)=q$, due to Proposition \ref{CAT0-prop}(ii),
one has that $\angle_q(\gamma(t),z_1)\geq \pi/2,$ $t\in (0,1].$ The
first variation formula (see Proposition \ref{CAT0-prop}(iii))
yields that $$0\geq \cos\angle_q(\gamma(t),z_1)=\lim_{s\to
0^+}\frac{d(q,z_1)-d(\gamma(s),z_1)}{s}.$$ Since $d$ is convex, the
function $s\mapsto \frac{d(q,z_1)-d(\gamma(s),z_1)}{s}$ is
non-increasing. Combining the latter two facts, it follows that
$$0\geq \frac{d(q,z_1)-d(\gamma(s),z_1)}{s},\ s\in (0,1].$$
In particular, $d(q,z_1)\leq d(\gamma(s),z_1)$ for every $s\in
(0,1],$ which concludes the first part of (\ref{propro}); the second
relation is proved similarly. Now, from the non-expansiveness of the
projection $P_{{\rm Im}(\gamma)}$ (see Proposition
\ref{CAT0-prop}(ii)) and relation (\ref{propro}), we obtain
$$d(q,\tilde q)=d( P_{{\rm Im}(\gamma)}(z_1),P_{{\rm
Im}(\gamma)}(z_2))\leq d( z_1,z_2).$$
\end{remark}

\noindent {\bf Proof of Theorem \ref{Riemann-double-proj}.} $"{\rm
(iii)\Rightarrow (i)\& (ii)}"$ If the Riemannian manifold $(M,g)$
has non-positive sectional curvature, $(M,d_g)$ is an Alexandrov NPC
space, see Bridson and Haefliger \cite[Theorem 1A.6]{Br-Ha}.
Consequently, by Theorem \ref{Ptolemaic-theorem}, $(M,d_g)$ has the
double-projection property. Moreover, by Proposition
\ref{CAT0-prop}(ii) it follows that the projective non-expansiveness
property also holds.

 $"{\rm
(i)\Rightarrow (iii)}"$ We assume that $(M,d_g)$ satisfies the
double-projection property, i.e.,  every $p\in M$ has a neighborhood
$U\subset M$ such that $(U,d)$ is a geodesic length space, and for
every two  geodesic convex, $U-$proximinal sets $S_1,S_2\subset U$,
the statements $(DP_1)$ and $(DP_2)$ are equivalent.

Let $p\in M$ be fixed and $B_g(p,\tilde \rho_p)\subset U$ be a
totally normal ball of $p,$ see do Carmo \cite[Theorem
3.7]{doCarmo}. Clearly, $B_g(p,\tilde \rho_p)$ inherits the above
properties of $U$. Fix also ${ W}_0,{ V}_0\in T{_pM\setminus
\{0\}}.$ We claim that the sectional curvature of the
two-dimensional subspace $\mathcal S=$ span$\{ { W}_0,{ V}_0
\}\subset T{ _pM}$ at ${ p}$ is non-positive. One may assume without
loss of generality that ${ V_0}$ and ${ W_0}$ are ${
g}-$perpendicular, i.e., ${ g}({ W}_0,{ V}_0)=0$.

Let $\kappa$ be an upper bound for the sectional curvature over the
closed ball $B_g[p,\tilde \rho_p]=\{q\in M:d_g(p,q)\leq \tilde
\rho_p\}$, and let $\kappa_1=\max\{1,\kappa\}.$  We fix $\delta>0$
such that
\begin{equation}\label{totally}
  \delta(\|{ W}_0\|_{ g}+\|{ V}_0\|_{ g})<\frac{1}{2}\min\left\{\tilde \rho_{ p},\frac{\pi}{\sqrt{\kappa_1}}\right\}.
\end{equation}

Let  $\sigma:[0,\delta]\to { M}$ be the geodesic segment
$\sigma(t)=\exp_{ p}(t{ V}_0)$ and ${ W}$ be the unique parallel
vector field along $\sigma$ with the initial data ${ W}(0)={ W}_0$.
For any $t\in [0,\delta]$, we define the geodesic segment
$\gamma_t:[0,\delta]\to { M}$ by
 $\gamma_t(u)=\exp_{\sigma(t)}(u{ W}(t)).$
Having in our mind these notations, we claim that
\begin{equation}\label{proj-proj}
    P_{{\rm Im}(\gamma_t)}(p)=\sigma(t)%\ {\rm  and}\ P_{{\rm Im}(\gamma_0)}(\sigma(t))=p
    \ {\rm for\ every}\ t\in [0,\delta].
\end{equation}
To show this, fix $t\in [0,\delta].$  Due to (\ref{totally}), for
every $u\in [0,\delta]$, the geodesic segment $\gamma_t|_{[0,u]}$
belongs to the normal ball
 $B{ _g}({ p},\tilde \rho_{ p})$;
thus, $\gamma_t|_{[0,u]}$ is the unique minimal geodesic joining the
point $\gamma_t(0)=\sigma(t)$ to $\gamma_t(u).$ Moreover, since ${
W}$ is the parallel transport of ${ W}(0)={ W}_0$ along $\sigma$, we
have ${ g}({ W}(t),\dot \sigma(t))={ g( W}(0),\dot \sigma(0)) ={
g}({ W}_0,{ V}_0)=0;$ therefore,
\begin{equation}\label{kilencven-fok}
    { g}(\dot \gamma_t(0),\dot \sigma(t)) ={ g( W}(t),\dot \sigma(t))=0.
\end{equation}
Since ${{\rm Im}(\gamma_t)}$ is compact, $P_{{\rm
Im}(\gamma_t)}(p)\neq \emptyset;$ let $q\in P_{{\rm
Im}(\gamma_t)}(p)$, and assume that $q\neq \sigma(t)$. It is clear
that the geodesic triangle $\Delta(p,q,\sigma(t))$ is included into
$B_g(p,\tilde \rho_p)$, and on account of (\ref{totally}), its
perimeter satisfies the inequality
\begin{equation}\label{perimeter}
d_g(p,q)+d_g(q,\sigma(t))+d_g(p,\sigma(t))<\frac{\pi}{\sqrt{\kappa_1}}.
\end{equation}
Moreover, due to the fact that $q\in P_{{\rm Im}(\gamma_t)}(p)$ and
(\ref{kilencven-fok}), the angles in the geodesic triangle
$\Delta(p,q,\sigma(t))$ fulfill
\begin{equation}\label{angles-geodesic}
  \measuredangle q\geq \frac{\pi}{2}\ {\rm and}\  \measuredangle \sigma(t)=\frac{\pi}{2}. \
\end{equation}
Now, we are in the position to apply Toponogov's comparison theorem
for triangles (where the curvature is bounded from above by the
number $\kappa_1>0$), see Klingenberg \cite[Theorem 2.7.6]{Klin}.
Namely, if $\overline \Delta(\overline p,\overline
q,\overline{\sigma(t)})$ is the comparison triangle for
$\Delta(p,q,\sigma(t))$ on the two-dimensional sphere with radius
$\frac{1}{\sqrt{\kappa_1}}$, the comparison angles in $\overline
\Delta(\overline p,\overline q,\overline{\sigma(t)})$ are not
smaller than their corresponding angles in $\Delta(p,q,\sigma(t))$.
Combining this fact with (\ref{angles-geodesic}), we get that
$$\measuredangle \overline q\geq \frac{\pi}{2}\ {\rm and}\   \measuredangle \overline{\sigma(t)}\geq \frac{\pi}{2}. \ $$
By the cosine rule for sides of a spherical triangle, the latter
inequalities yield
$$\cos d_g(p,\sigma(t))-\cos d_g(p,q)\cos d_g(q,\sigma(t))=\sin d_g(p,q)\sin d_g(q,\sigma(t))\cos \overline q \leq 0;$$
$$\cos d_g(p,q)-\cos d_g(p,\sigma(t))\cos d_g(q,\sigma(t))=\sin  d_g(p,\sigma(t))\sin d_g(q,\sigma(t))\cos \measuredangle \overline{\sigma(t)} \leq 0.$$
Adding these inequalities and rearranging them, we obtain
$$\left[1-\cos d_g(q,\sigma(t))\right]\cdot[\cos d_g(p,q)+\cos
d_g(p,\sigma(t))]\leq 0,$$ which is equivalent to
$$\sin^2\frac{d_g(q,\sigma(t))}{2}\cos\frac{d_g(p,q)+d_g(p,\sigma(t))}{2}\cos\frac{d_g(p,q)-d_g(p,\sigma(t))}{2}\leq 0.$$
Since $q\neq \sigma (t)$, the first term is positive. On account of
(\ref{perimeter}), the third term is also positive. Thus, the second
term is necessarily non-positive, i.e.,
$d_g(p,q)+d_g(p,\sigma(t)))\geq \pi,$ which contradicts
(\ref{perimeter}). Consequently, $P_{{\rm Im}(\gamma_t)}(p)$
contains the unique element $\sigma(t)$, which concludes the proof
of  (\ref{proj-proj}).

In the same way as in (\ref{proj-proj}), we can prove
\begin{equation}\label{proj-proj-2}
     P_{{\rm Im}(\gamma_0)}(\sigma(t))=p
    \ {\rm for\ every}\ t\in [0,\delta].
\end{equation}
Thus, we can conclude from (\ref{proj-proj}) and (\ref{proj-proj-2})
that for every $t\in [0,\delta]$, $$P_{{\rm Im}(\gamma_0)}(P_{{\rm
Im}(\gamma_t)}(p))=p,$$ i.e., $(DP_1)$ holds for the point  $p\in
{\rm Im}(\gamma_0)$ and sets $S_1={\rm Im}(\gamma_0)$ and $S_2={\rm
Im}(\gamma_t)$, respectively.  Since these sets are geodesic convex
and compact (thus, $B_g(p,\tilde \rho_p)-$proximinal), the validity
of the double-projection property implies that
 $(DP_2)$ holds  too, i.e., $p$ is the best approximation point from Im$(\gamma_0)$ to
 Im$(\gamma_t)$. Formally, we have
 $$d_g(p,P_{{\rm Im}(\gamma_t)}(p))\leq d_g(z_1,z_2)\ {\rm for\  all}\ (z_1,z_2)\in
{\rm Im}(\gamma_0)\times {\rm Im}(\gamma_t)\  {\rm and} \ t\in
[0,\delta].$$ In particular, for every $t,u\in [0,\delta]$, we have
\begin{equation}\label{egyenlotelneseg-utolso}
    d_g(p,\sigma(t))
 \leq d_g( \gamma_0(u), \gamma_t(u)).
\end{equation}
By using  the parallelogramoid of Levi-Civita for calculating the
sectional curvature $K_p(\mathcal S)$ at $p$ and for the
two-dimensional subspace $\mathcal S$=span$\{ { W}_0,{ V}_0
\}\subset T_pM$, see Cartan \cite[p. 244-245]{Cartan}, we obtain
from (\ref{egyenlotelneseg-utolso}) that
$$
  K_p(\mathcal S)=\lim_{u,t\to 0}\frac{d_g^2(p,\sigma(t))-d_g^2(
\gamma_0(u), \gamma_t(u))}{d_g( p, \gamma_0(u))\cdot d_g( p,
\sigma(t))}\leq 0.$$ This concludes the proof of  $"{\rm
(i)\Rightarrow (iii)}"$.

 $"{\rm
(ii)\Rightarrow (iii)}"$ Let us keep the notations and constructions
from above. A similar geometric reasoning as in the proof of
(\ref{proj-proj}) yields that
\begin{equation}\label{utolosoookkok}
 P_{{\rm Im}(\sigma)}(\gamma_t(u))=\sigma(t)
    \ {\rm for\ every}\ t,u\in [0,\delta].
\end{equation}
Since $S={\rm Im}(\sigma)$ is a geodesic convex $B_g(p,\tilde
\rho_p)-$proximinal set and the projection non-expansiveness
property holds,  on account of (\ref{non-exp}) and
(\ref{utolosoookkok}) we obtain for every $t,u\in [0,\delta]$ that
$$d_g(p,\sigma(t))=d_g(\sigma(0),\sigma(t))=d_g(P_{{\rm Im}(\sigma)}(\gamma_0(u)),P_{{\rm Im}(\sigma)}(\gamma_t(u)))\leq d_g( \gamma_0(u), \gamma_t(u)),$$
which is nothing but relation (\ref{egyenlotelneseg-utolso}). It
remains to follow the previous proof.
  \hfill $\diamondsuit$

\section{Proof of Theorem \ref{Minkowski-theorem} and the double-projection property on non-reversible Finsler-Minkowski spaces}\label{sect-4}

\noindent  {\bf Proof of Theorem \ref{Minkowski-theorem}.} Let
$S_1,S_2\subset \mathbb R^n$ be two convex and $\mathbb
R^n-$proximinal sets. Note that the implication  $"(DP_2)\Rightarrow
(DP_1)"$ is proved analogously as in Theorem
\ref{Ptolemaic-theorem}.

 Let us prove  $"(DP_1)\Rightarrow (DP_2)"$. To do this, let $q\in
 S_1$ such that $q\in P_{S_1}(P_{S_2}(q)).$ Due to Proposition
 \ref{Minkowski-prop}(i), both sets $S_1$ and $S_2$ are $\mathbb
 R^n-$Chebishev. Consequently, there exists a unique element $\tilde q\in
 S_2$ such that $P_{S_2}(q)=\tilde q$ and $P_{S_1}(\tilde q)=q.$ On
account of Proposition
 \ref{Minkowski-prop}(ii), the latter relations  are equivalent to  $$g_{\tilde
q -q}(\tilde q -q,z_1-q)\leq 0 \ {\rm for\ all}\ z_1\in S_1;$$
$$g_{\tilde q -q}( q -\tilde q,z_2-\tilde q )\leq 0 \ {\rm for\ all}\
z_2\in S_2.$$ Adding these inequalities, we obtain $g_{\tilde q
-q}(\tilde q -q,\tilde q -q-z_2+z_1)\leq 0.$ By applying the
fundamental inequality (see Proposition
 \ref{Minkowski-prop}(iii)) and relation (\ref{F-metrika}), we have
\begin{eqnarray*}
d_F^2(q,\tilde q )&=&F^2(\tilde q -q)=g_{\tilde q -q}(\tilde q
-q,\tilde q -q)\\&\leq& g_{\tilde q -q}(\tilde q -q,z_2-z_1)\\&\leq&
F(\tilde q -q)\cdot F(z_2-z_1)\\&=&d_F(q,\tilde q )\cdot
d_F(z_1,z_2).
\end{eqnarray*}
Consequently, $d_F(q,\tilde q )\leq d_F(z_1,z_2)$ for every $z_1\in
S_1$ and $z_2\in S_2$, which means that $q\in S_1$ is the best
approximation element from  $S_1$ to $S_2$.
 $\hfill $ $\diamondsuit$

\begin{remark}\rm

Let $(\mathbb R^n,F)$ be a not necessarily reversible
Finsler-Minkowski space; the metric distance $d_F$ is usually only a
quasi-metric. Even in this case, it is possible to state a similar
result as Theorem \ref{Minkowski-theorem}, slightly reformulating
the double-projection property.

 Let
$S_1,S_2\subset \mathbb R^n$ be two convex and $\mathbb
R^n-$proximinal sets, and $q\in S_1$. Note that $S_1,S_2$ are
forward and backward $\mathbb R^n-$Chebishev sets, see Proposition
\ref{Minkowski-prop}(i). In the {\it forward} case, we consider the
following statements:
\begin{itemize}
\item[$(DP_1^+)$]  $q= (P_{S_1}^-\circ
P_{S_2}^+)(q)$;
\item[$(DP_2^+)$]   $d_F(q,P_{S_2}^+(q))\leq d_F(z_1,z_2)\   {\rm for\ all}\
z_1\in S_1,\ z_2\in S_2.$
\end{itemize}
In the {\it backward} case, we
 consider similar statements:
\begin{itemize}
\item[$(DP_1^-)$]  $q= (P_{S_1}^+\circ P_{S_2}^-)(q)$;
\item[$(DP_2^-)$]  $d_F(P_{S_2}^-(q), q)\leq d_F(z_2,z_1)\  {\rm for\ all}\
z_1\in S_1,\ z_2\in S_2$.
\end{itemize}
\end{remark}

Exploiting Proposition
 \ref{Minkowski-prop}(ii) in its full generality,  we can  show as  in Theorem
 \ref{Minkowski-theorem}:
\begin{theorem}\label{Minko-67}
Let $(\mathbb R^m,F)$ be a Finsler-Minkowski space. Then for every
two convex and $\mathbb R^n-$proximinal sets $S_1,S_2\subset \mathbb
R^n$, we have:
\begin{itemize}
  \item[{\rm (i)}] $(DP_1^+)\Leftrightarrow (DP_2^+);$
  \item[{\rm (ii)}] $(DP_1^-)\Leftrightarrow (DP_2^-).$
\end{itemize}
\end{theorem}

\begin{remark}\rm
Usually, the map $P_{S_1}^-\circ P_{S_2}^+$ in $(DP_1^+)$ {\it
cannot} be replaced  either by $P_{S_1}^+\circ P_{S_2}^+$ or by
$P_{S_1}^+\circ P_{S_2}^-$ or by $P_{S_1}^-\circ P_{S_2}^-$. (The
same is true for $P_{S_1}^+\circ P_{S_2}^-$ in  $(DP_1^-)$.) In
order to give a concrete example, we recall the Matsumoto norm, see
\cite{Matsumoto}, which describes the walking-law on a mountain
slope (under the action of gravity), having an angle $\alpha\in
[0,\pi/2)$ with the horizontal plane. The explicit form of this norm
$F:\mathbb R^2\to [0,\infty)$ is
\begin{equation}\label{Matsumoto_metrika} F(y)=\left\{
  \begin{array}{ll}
    \frac{y_{1}^{2}+y_{2}^{2}}{v\sqrt{y_{1}^{2}+y_{2}^{2}%
}+\frac{g}{2}y_{1}\sin\alpha}, &   y=(y_{1},y_{2})\in\mathbb{R}^{2}\setminus\{(0,0)\}; \\
    0, &  y=(y_{1},y_{2})=(0,0),
  \end{array}
\right.
\end{equation}
where $v$ [$m/s$] is the constant speed on the horizontal plane,
$g\approx9.81$ $[{m}/{s^2}],$ and $g\sin\alpha \leq v$. The pair
$(\mathbb R^2,F)$ is a typical non-reversible Finsler-Minkowski
space, and it becomes reversible if and only if $\alpha=0$.

\begin{figure*}
\begin{center}
%\begin{center}
 \includegraphics[width=0.6\textwidth]{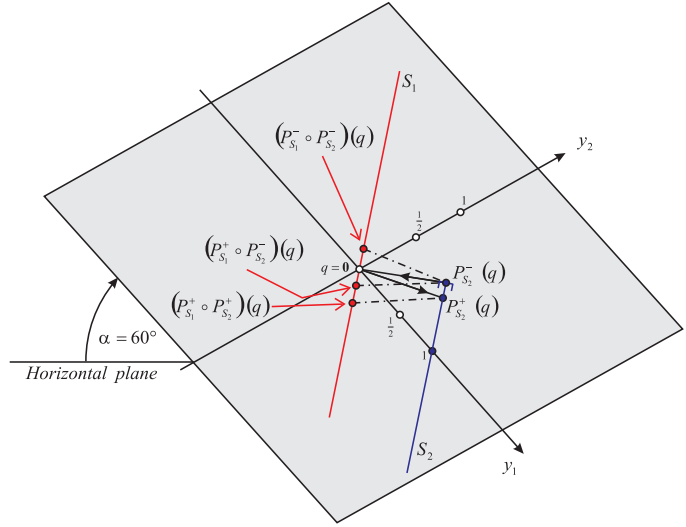}
\quote{\caption{Apart from the case $P_{S_1}^-\circ P_{S_2}^+$, the
compositions of forward and/or backward metric projections at
$q=(0,0)$ are scattered away from  $q$.}}
\label{fig:ama}       % Give a unique label
\end{center}
\end{figure*}

 Let $v=10$ and $\alpha=\pi/3$ in  (\ref{Matsumoto_metrika}), and
 consider the convex and closed sets
$$\left.
\begin{array}{lcl}
S_1=\{(y_1,y_2)\in \mathbb R^2: y_1+y_2=0\},\\ S_2=\{(y_1,y_2)\in
\mathbb R^2: y_1+y_2=1,\ y_1\geq 1/2\}.
\end{array}\right.$$
  Let also $q=(0,0)$. A direct
calculation yields  $(P_{S_1}^-\circ P_{S_2}^+)(q)=q$, and
$d_F(q,P_{S_2}^+(q))\leq d_F(z_1,z_2)\   {\rm for\ all}\  z_1\in
S_1,\ z_2\in S_2.$ However,
 we have

$(P_{S_1}^+\circ P_{S_2}^+)(q)=(0.32338512,-0.32338512)\neq q,$

$(P_{S_1}^+\circ
P_{S_2}^-)(q)=P_{S_1}^+(1/2,1/2)=(0.23349577,-0.23349577)\neq q,$

$(P_{S_1}^-\circ
P_{S_2}^-)(q)=P_{S_1}^-(1/2,1/2)=(-0.08988935,0.08988935)\neq q,$

\noindent see also Figure 1.
\end{remark}

\vfill\eject

\vspace{0.5cm} \noindent {\bf Acknowledgment.} A. Krist\'aly was
supported by the grant CNCS-UEFISCDI/PN-II-RU-TE-2011-3-0047 and
%by the Romanian National Authority for Scientific
%Research, CNCS-UEFISCDI, project number PN-II-ID-PCE-2011-3-0241 and
by the J\'anos Bolyai Research Scholarship. D. Repov\v s was
supported by a grant of the Slovenian Research Agency P1-0292-0101.

\end{document}